\documentclass[12pt,reqno]{article}

\usepackage[usenames]{color}
\usepackage{amssymb}
\usepackage{graphicx}
\usepackage{amscd}

\usepackage[colorlinks=true,
linkcolor=webgreen,
filecolor=webbrown,
citecolor=webgreen]{hyperref}

\definecolor{webgreen}{rgb}{0,.5,0}
\definecolor{webbrown}{rgb}{.6,0,0}

\usepackage{color}
\usepackage{fullpage}
\usepackage{float}

\usepackage{graphics,amsmath,amssymb}
\usepackage{amsthm}
\usepackage{amsfonts}
\usepackage{latexsym}
\usepackage{epsf}

\usepackage{rotating}

\newcommand{\bsq}{\vrule height .9ex width .8ex depth -.1ex}

\newcommand{\sE}{{\cal E}}

\newcommand{\sI}{{\cal I}}

\newcommand{\sM}{{\cal M}}

\newcommand{\eeq}{\end{equation}}

\newcommand{\ra}{\rightarrow}

\newcommand{\beql}[1]{\begin{equation}\label{#1}}
\newcommand{\eqn}[1]{(\ref{#1})}
\newcommand{\beq}{\begin{displaymath}}
\newcommand{\eeqno}{\end{displaymath}}

\newcommand{\qandq}{\quad\mbox{and}\quad}

\newcommand{\qasq}{\quad\mbox{as}\quad}

\newcommand{\seqnum}[1]{\href{http://www.research.att.com/cgi-bin/access.cgi/as/~njas/sequences/eisA.cgi?Anum=#1}{\underline{#1}}}

\begin{document}

\begin{center}
\epsfxsize=4in
\end{center}

\begin{center}
\vskip 1cm{\LARGE\bf Integer Sequences from Queueing Theory}
\vskip 1cm \large
Joseph Abate\\
900 Hammond Road\\
Ridgewood, NJ 07450-2908\\
USA\\
\ \\
Ward Whitt\\
Department of Industrial Engineering and Operations Research\\
Columbia University\\
New York, NY 10027-6699\\
USA\\
\href{mailto:ww2040@columbia.edu}{\tt ww2040@columbia.edu} \\
\end{center}

\vskip .2 in

\begin{abstract}
  Operators on probability distributions can be expressed
  as operators on the associated moment sequences, and so
correspond to operators on integer sequences. Thus, there is an
opportunity to apply each theory to the other.
  Moreover, probability models can be sources
of integer sequences, both classical and new, as we show by
considering the classical $M/G/1$ single-server queueing model. We
identify moment sequences that are integer sequences.  We
establish connections between the $M/M/1$ busy period distribution and the Catalan and Schroeder numbers.
\end{abstract}

\newtheorem{theorem}{Theorem}
\newtheorem{corollary}[theorem]{Corollary}
\newtheorem{lemma}[theorem]{Lemma}
\newtheorem{proposition}[theorem]{Proposition}
\newtheorem{conjecture}[theorem]{Conjecture}
\newtheorem{defin}[theorem]{Definition}
\newenvironment{definition}{\begin{defin}\normalfont\quad}{\end{defin}}
\newtheorem{examp}[theorem]{Example}
\newenvironment{example}{\begin{examp}\normalfont\quad}{\end{examp}}
\newtheorem{rema}[theorem]{Remark}
\newenvironment{remark}{\begin{rema}\normalfont\quad}{\end{rema}}

\providecommand{\e}[1]{\ensuremath{\times 10^{#1}}}

\section{Introduction}\label{secIntro}
Given the well established link between integer sequences and combinatorics, it is surprising that
there are so few connections between integer sequences and queueing theory,
because two prominent scholars, John Riordan (1902--1988) and Lajos Tak\'{a}cs (1924-- \, ), were active in both
combinatorics and queueing theory \cite{R58,R62,R68,T62,T67,T90}.
Of course, these authors saw connections; e.g.,
the index of Riordan's queueing book \cite{R62} contains entries for Bell polynomials, binomial moments,
Catalan numbers, cumulant generating function, Lagrange expansion, and rooted trees.

In the OEIS \cite{S10}, a search on ``Riordan'' yields 1221 entries, while a search on ``queueing'' produces 2.
We think Riordan would want less disparity.
One exception is the entry submitted by A. Harel (A122525) making connection to the Erlang delay formula associated with the classical $M/M/s$ queue.
A search on ``queue'' produces 15 entries, but these primarily focus on combinatorial problems \cite{A04,S05}
or a queue as a tool in computer science rather than queueing theory \cite{AB94};
queueing theory concerns stochastic process models describing congestion, e.g., the probability distribution of customer waiting times \cite{R62,T62}.

The purpose of this paper is to point out connections between the theories of probability and integer sequences,
and to exhibit some integer sequences that arise in queueing theory.
As a branch of probability theory, queueing theory has exploited
the classical analytical approach to probability theory
using transforms \cite{R62,T62}, whose series expansions involve the integer moments of the probability distributions
and close relatives.  The entire probability distribution
is characterized by the sequence of moments in great generality \cite{A65,DS97}.
Thus we were motivated to introduce an operational calculus for manipulating probability distributions on the positive halfline by either manipulating
the associated Laplace transforms or
the associated moment sequences \cite{AW96}; a quick overview is provided by \cite[Tables 1-3]{AW96}.
A more recent paper in the same spirit is \cite{FG00}.

The operational calculus for probability distributions on the positive halfline in the framework of moment sequences is closely related to operators
commonly used to analyze integer sequences.
Since many operators on integer
sequences can be applied to moment sequences arising in
probability theory, both when they are integers and when they are not, there is an opportunity for experts on
integer sequences to contribute to probability theory through moment sequences.
(The connection is also discussed in \cite{B09}.)  The
probability connection also provides concrete models where integer
sequences arise.

Here is how the present paper is organized.  We start in \S
\ref{secMoment} by reviewing moment sequences in probability
theory and relating operators on probability distributions to
operators on sequences.
Then in \S \ref{secMM1} we review the
classical $M/G/1$ single-server queueing model and identify
integer sequences arising in that context. We consider the $M/M/1$
busy period distribution, its stationary excess distribution and the equilibrium time to emptiness. In
that context we identify random quantities whose probability
density functions have the Catalan and Large Schroeder numbers as
moments. Finally, we consider the moment sequence for the $M/G/1$
steady-state waiting time. We give the proofs of Theorems \ref{th1} and \ref{th2} in
\S \ref{secProofs}.

\section{Moment Sequences in Probability}\label{secMoment}

  Let $Z$ be a nonnegative random variable with {\em cumulative distribution function} (cdf) $F$
and {\em probability density function} (pdf) $f$, i.e.,
\beq
F(t) \equiv P(Z \le t) = \int_{0}^{t} f(u) \, du, \quad t \ge 0,
\eeqno
where $\equiv$ means ``defined as.''
Let $\hat{f} (s)$ be the {\em Laplace transform} (LT) of $f$ (and thus $Z$) and let $\phi (x)$ be the associated {\em moment generating function} (mgf)
of $f$, i.e.,
\begin{eqnarray}\label{a2}
\hat{f}(s) & \equiv & \int_{0}^{\infty} e^{-st} f(t) \, dt \equiv E[e^{-sZ}] \qandq\nonumber \\
\phi (x)   & \equiv & \hat{f} (-x) = \int_{0}^{\infty} e^{x t} f(t) \, dt \equiv E[e^{xZ}],
\end{eqnarray}
where $s$ in \eqn{a2} is understood to be a complex number with positive real part, while $x$ in \eqn{a2}
is a positive real number.  The LT is always well defined; we assume that there exists $x^* > 0$ such that
$\phi (x) < \infty$ for $x < x^*$.

We obtain sequences by considering series expansions for the mgf $\phi$; in particular, we can write
\beql{a3}
\phi(x) = \sum_{n = 0}^{\infty} \frac{m_n x^n}{n!} = \sum_{n = 0}^{\infty} \mu_n x^n,
\eeq
where of course we must have $\mu_n = m_n/n!$, $n \ge 0$.
From the probability perspective, 
the object of primary interest is the cdf $F$, but the moments $m_n$ provide a useful partial characterization.
The moment sequence can be calculated and employed to derive other quantities of interest in probability models
\cite{ACLW95,AW96,CL96}.  There is a developing operational calculus for manipulating probability distributions via their moment sequences \cite{AW96}.

From the sequence perspective, we may instead regard the sequences
$\{m_n: n \ge 1\}$ and $\{\mu_n: n \ge 1\}$ as the objects of
primary interest.  The two perspectives meet with the mgf $\phi$.
From the sequence perspective, $\phi (x)$ arises as the {\em
generating function} (gf) of the sequence $\{\mu_n: n \ge 1\}$,
and we speak of $\{\mu_n\}$ as being its {\em coefficients}, while
$\phi (x)$ arises as the {\em exponential generating function}
(egf) of the sequence $\{m_n: n \ge 1\}$.  Of course, $\{m_n: n
\ge 1\}$ is always an integer sequence whenever $\{\mu_n: n \ge
1\}$ is, but not necessarily conversely. We will be giving
examples where both are integer sequences.

A simple example from probability theory is the exponential distribution with mean $m$, specified by
\beq
F(t) \equiv 1 - e^{-t/m}, \quad f(t) \equiv (1/m) e^{-t/m}, \quad t \ge 0, \qandq \phi (x) \equiv (1 - m x)^{-1},
\eeqno
which has associated sequences
\beq
m_n \equiv n! m^n \qandq \mu_n \equiv m^n, \quad n \ge 0.
\eeqno
It is immediate that $\{\mu_n\}$ and $\{m_n\}$ are both elementary integer sequences whenever $m$ is an integer.
Since the mean (first moment) is probabilistically only a scale parameter, depending on how the measurement units are defined,
it is natural to follow the convention that the first moment is $1$; here that gives $m_1 = 1$.  Then we obtain the fundamental integer sequences
$m_n = n!$ and $\mu_n = 1$, $n \ge 1$.

As a consequence of the relations outlined above, we see that results about probability distributions can be translated into
results about integer sequences, provided that the sequence $\{\mu_n\}$ or $\{m_n\}$ is indeed an integer sequence.  Similarly,
results about integer sequences can be translated into
results about probability distributions, provided that the egf or gf of the integer sequence is indeed the mgf of a
bonafide pdf.

We first observe that there is a natural probabilistic setting in which, not only is $\{m_n\}$ a moment sequence, but so is the
associated sequence $\{\mu_n\}$.  That occurs in spectral representations, where one pdf serves as a mixing pdf for another.
In particular, suppose that the pdf $f$ can be represented as a continuous mixture of exponential pdf's via
\beql{a31}
f(t) = \int_{\tau_1}^{\tau_2} y^{-1} e^{-t/y} g(y) \, dy, \quad t \ge 0,
\eeq
in which case we call $g$ the {\em mixing pdf} for $f$ \cite{AW88b}.
Append a superscript $f$ to $\{m_n\}$ and $\{\mu_n\}$ to denote dependence upon $f$.
We observe that the sequence $\{\mu_n\}$ is itself the moment sequence of the associated mixing pdf $g$.  Hence we call $\{\mu_n\}$ the {\em mixing moments} of $f$.
(We omit the elementary proofs in this section.)
\begin{proposition}{$($mixing moments$)$}\label{propMixing} For pdf's $f$ and $g$ related via {\em \eqn{a31}},
$m^g_n = \mu^f_n$, $n \ge 1$.
\end{proposition}

A canonical operation in probability theory is convolution.  If
two independent nonnegative random variables $Z_1$ and $Z_2$ with
pdf's $f_{Z{_1}}$ and $f_{Z{_2}}$ are added, then the sum $Z_1 +
Z_2$ has a pdf that is the {\em convolution} of the two component
pdf's, i.e.,
\beq
 f_{Z{_1} + Z{_2}} (t) = \int_{0}^{t}
f_{Z{_1}} (y) f_{Z{_2}} (t - y) \, dy, \quad t \ge 0.
\eeqno
One reason
transforms are so frequently used in probability is that they
convert convolution into a simple product; i.e., the associated
mgf's are related by
\beq
 \phi_{Z{_1} + Z{_2}} (x) =
\phi_{Z{_1}} (x) \phi_{Z{_2}} (x).
\eeqno
Because of the assumed stochastic independence, the moments are related by
the binomial theorem,
\beq
 m^{Z_{_1} + Z{_2}}_n =
\sum_{k=1}^{n} \binom{n}{k} m^{Z{_1}}_k m^{Z_{_2}}_{n-k}.
\eeqno
Thus, if $\{ m^{Z{_1}}_n\}$ and $\{ m^{Z{_2}}_n\}$ are integer
sequences, then so is $\{m^{Z{_1} + Z{_2}}_n\}$. For example, by
above, that occurs when $Z_1$ and $Z_2$ have exponential
or deterministic distributions with integer
means.

With this background, we can interpret \cite[Tables 1-3]{AW96}, which show various operators mapping one probability distribution
into another.  There are four columns:  The first column contains the name and notation for the operator; the second column shows how the operator acts on
LT's; the third column shows how it acts on pdf's; and the fourth column shows how it acts on moment sequences.
From the perspective of integer sequences, the fourth column shows how it acts on the coefficients of the egf;
we could then add a fifth column showing how it acts
on the corresponding coefficients of the gf.
Those familiar with integer sequences might want to translate the LT into the mgf by replacing $s$ with $-x$, and then interpret
that mgf as the egf of the given $k^{\rm th}$ moment.

In the rest of this section we highlight a few striking connections between operators on probability distributions, as in \cite{AW96},
and operators on integer sequences.
First, a standard operator on integer sequences is the simple {\em shift to the left}, e.g., converting
$1, 1, 2, 6, 22, 90, \ldots$ to $1, 2, 6, 22, 90, \ldots$.
  We now show that,
probabilistically, the simple shift applied to the coefficients $\mu_n$ corresponds to constructing the stationary-excess cdf
of a cdf on the positive halfline having mean $1$.

Given a nonnegative real-valued random variable $Z$ with cdf $F$ having
finite moments $m_k$, $k \ge 1$, let $Z_e$ be a random variable with the associated {\em stationary-excess cdf} $F_e$
(a.k.a. the equilibrium excess or stationary residual-life cdf),
defined by
\beql{a21}
F_e (t) \equiv P(Z_e \le t) \equiv \frac{1}{m_1}\int_{0}^{t} (1 - F(u)) \, du, \quad t \ge 0.
\eeq
The stationary-excess cdf frequently arises in renewal theory; see  of
\cite[Examples 7.16,7.17, 7.23, 7.24]{R07} and \cite{W85}; it appears in \cite[Table 2]{AW96}.
(A search on ``renewal theory'' in the OEIS gives two unrelated entries.)
For us, the important fact is that the
random variable $Z_e$ has moments
\beql{a22}
m_{e,k} \equiv E[Z_e^k] = \frac{m_{k+1}}{(k+1) m_1}, \quad k \ge 1.
\eeq
Hence, the transformation from a cdf of a nonnegative random variable to its associated
stationary-excess cdf produces a simple shift on the gf coefficients.
\begin{proposition}{$($the stationary-excess operator$)$}\label{propShift} Let $F$ be the cdf of a
nonnegative random variable with mean $1$, mgf $\phi$ in {\em \eqn{a2}} and
associated sequence of mixing moments $\{\mu_n\}$ in {\em \eqn{a3}} $($coefficients of $\phi (x)$ when it is regarded as a gf$)$.
The associated stationary-excess cdf $F_e$ in {\em \eqn{a21}} has mgf
$\phi_{e} (x) = (\phi (x) - 1)/x$.
The mixing moments of $\phi_e (x)$ $($coefficients of $\phi_e$ when it is regarded as a gf$)$ are
$\mu_{e,k} = \mu_{k+1}$, $k \ge 1$.
\end{proposition}

A similar relationship holds for the {\em stationary-lifetime operator}, mapping a pdf $f$ into the associated pdf
\beql{a25}
f_s (t) \equiv \frac{t f(t)}{m_1}, \quad t \ge 0.
\eeq
The stationary-lifetime pdf also frequently arises in renewal theory, e.g., \cite[\S 7.7]{R07}, and also appears in \cite[Table 2]{AW96}.
For us, the important fact is that the
moments are related by
\beq
m_{s,k} \equiv \frac{m_{k+1}}{m_1}, \quad k \ge 1,
\eeqno
just like \eqn{a22} without the $k+1$ in the denominator.
Hence, the transformation from a pdf of a nonnegative random variable to its associated
stationary-lifetime pdf produces a simple shift on the moments $m_n$ (the egf coefficients).
\begin{proposition}{$($the stationary-lifetime operator$)$}\label{propShift2}
Let $f$ be the pdf of a nonnegative random variable with mean $1$, mgf $\phi$ in {\em \eqn{a2}} and
associated sequence of moments $\{m_n\}$ in {\em \eqn{a3}}.  The associated stationary-lifetime pdf $f_s$ in {\em \eqn{a25}} has mgf
$\phi_{s} (x) = \phi' (x)$.
The coefficients of $\phi_s$ when it is regarded as an egf are
$m_{s,k} = m_{k+1}$, $k \ge 1$.
\end{proposition}

We now turn to another basic probability operator, which is conveniently related to a
continued fraction representation of the mgf \cite{AW99a,AW99b, B09, CBVWJ08}. 
If $\hat{f}$ is the LT of a pdf $f$, then the associated {\em exponential mixture pdf} has LT
\beql{a41}
\hat{f}_{\sE\sM} (s) \equiv (1 + s \hat{f}(s))^{-1};
\eeq
it appears in \cite[Table 3]{AW96}.
The special case of exponential mixtures of inverse Gaussian (EMIG) distributions is discussed in \cite[\S 8]{AW96} and in \cite{AW99b}.

The probabilistic exponential mixing operator has a simple manifestation in the continued fraction representation of the mgf,
regarding that mgf as a gf.
Starting with the LT $\hat{f} (s)$ of a pdf $f$, if we represent the associated mgf as a formal power series
by
\beql{power}
\phi (x) \equiv \hat{f} (-x) = 1 + \mu_1 x + \mu_2 x^2 + \mu_3 x^3 + \mu_4 x^4 + \mu_5 x^5 + \ldots ,
\eeq
then the {\em corresponding continued fraction} (CF) is
\beql{CF}
\phi (x) = \frac{1}{1-} \frac{h_1 x}{1-} \frac{h_2 x}{1-} \frac{h_3 x}{1-} \frac{h_4 x}{1-} \ldots
\eeq
where $h_1 \equiv \mu_1$, $h_2 \equiv (\mu_2 - \mu_1^2)/\mu_1$, etc.  When $h_n > 0$ for all $n$
we have an $S$-fraction and the underlying pdf $f$ is {\em completely monotone} (CM) \cite{AW99a}.

\begin{proposition}{$($the exponential-mixture operator$)$}\label{propShift3} Let $f$ be the pdf of a nonnegative random variable with mean $1$,
mgf $\phi$ in {\em \eqn{a2}} with associated sequence of CF coefficients $\{h_n\}$ in {\em \eqn{CF}}.
The associated exponential-mixture pdf $f_{\sE\sM}$ with LT in {\em \eqn{a41}} has CF coefficients
\beq
h_{\sE\sM,1} = 1, \quad h_{\sE\sM,n} = h_{n-1}, \quad n \ge 2;
\eeqno
i.e., it produces a shift to the right.  The inverse exponential mixture operator
{\em \cite[(7.3), p. \ 94]{AW96}} gives the corresponding shift to the left.
\end{proposition}

\paragraph{Proof.}  From the transform expression in \eqn{a41}, the conclusion is immediate:  Given the CF
representation for $\hat{f}(s)$, it is immediate that the corresponding CF for $(1 + s \hat{f}(s))^{-1}$
shifts the coefficients one to the right; i.e., we write the CF for $(1 + s \hat{f}(s))^{-1}$ as
$1/(1+s \hat{f} (s))$, inserting the CF for $\hat{f} (s)$.~~~\bsq


\section{Queueing Examples}\label{secMM1}

\subsection{The $M/G/1$ Model}

The $M/G/1$ queue is a basic model in queueing theory, usually discussed in queueing textbooks; e.g., \cite[\S 8.5]{R07}.  There is a single server
with unlimited waiting room.   Customers arrive according to a Poisson process (the first $M$, for Markov) with rate $\lambda$, $0 < \lambda < \infty$.
If the system is empty, then the customer goes immediately into service; otherwise the customer waits in queue.
The successive service times come from a sequence of independent and identically distributed (i.i.d.) random variables with cdf $G$
having mean $1/\mu$, $0 < \mu < \infty$.  We will mostly consider the easiest case, in which the service-time cdf $G$ is exponential;
then the model is denoted by $M/M/1$.

A waiting customer enters service immediately upon service completion.
Let $Q(t)$ be the number of customers in the system at time $t$ for $t \ge 0$.  In the $M/M/1$ model, the stochastic process $Q \equiv \{Q(t): t \ge 0\}$ is
a birth-and-death stochastic process, with constant birth rate $\lambda$ and constant death rate $\mu$.  Let $\rho \equiv \lambda/\mu$ be the {\em traffic intensity}.
If $\rho < 1$, then $P(Q(t) = j|Q(0) = i) \ra (1 - \rho)\rho^j$ as $t \ra \infty$ for each $i$ and $j$; i.e.,
$Q(t)$ converges in distribution to a geometric distribution on the nonnegative integers, having mean $\rho/(1- \rho)$.
We assume that $\rho < 1$, under which the system is said to be stable.
(If $\rho \ge 1$, then $P(Q(t) \le j|Q(0) = i) \ra 0$ as $t \ra \infty$ for each $i$ and $j$.)
For the $M/G/1$ model, the steady-state distribution of $Q(t)$ is characterized by the Pollaczek-Khintchine transform \cite{R62,T62}.

\subsection{The $M/M/1$ Busy Period}

A busy period is the time from the arrival of a customer finding an empty system until the system is empty again \cite[\S 4.8]{R62}.  The first passage time of $Q$
from any state $j> 0$ to $j - 1$ is distributed as a busy period.
Without loss of generality, we can measure time in units of mean service times, so that we let $\mu \equiv 1$.  Then the model has the single parameter $\rho \, (=\lambda)$.
Let $X_{\rho}$ denote a random variable with the busy period distribution, as a function of the traffic intensity $\rho$.  Let $B_{\rho}$ be the
cdf of $X_{\rho}$, i.e., $B_{\rho} (t) \equiv P(X_{\rho} \le t)$, $t \ge 0$,
and let $b_{\rho}$ be the associated pdf.  It turns out that
\beq
b_{\rho} (t) \equiv \frac{1}{t\sqrt{\rho}}e^{-(1 + \rho)t} I_{1}(2t\sqrt{\rho}), \quad t \ge 0,
\eeqno
where $I_{1} (t)$ is the Bessel function of the first kind \cite[(39), p. 63]{R62}.
The LT of $b_{\rho}$ (and thus $X_{\rho}$) is
\beq
\hat{b}_{\rho} (s) \equiv \int_{0}^{\infty} e^{-st} b_{\rho} (t) \, dt \equiv E[e^{-sX_{\rho}}] = \frac{1 + \rho + s - \sqrt{(1 + \rho +s)^2 - 4 \rho}}{2\rho},
\eeqno
\cite[(38), p. 63]{R62}.  As usual, the moments can be obtained by differentiating the Laplace transform.  The first two moments are
$E[X_{\rho}] = (1 - \rho)^{-1}$ and $E[X_{\rho}^2] = 2(1 - \rho)^{-3}$.

\subsection{The Busy-Period Moment Sequence After Scaling}

Our goal here is to obtain interesting integer sequences from the sequence of successive integer moments of a busy period $X_{\rho}$.
To obtain integer sequences, we first perform a change of variables,
introducing $\sigma \equiv \rho/(1 - \rho)$ (the mean steady state number in system) or, equivalently, $\rho \equiv \sigma/(1 + \sigma)$.
Notice that the first two moments become $1 + \sigma$ and $2 (1 + \sigma)^3$.  It turns out that the entire moment sequence becomes an integer sequence
whenever $\sigma$ is an integer.  Hence, from this first step, we obtain an entire sequence of integer sequences.

To seek simple integer sequences, we further scale these moment sequences all to have mean (first moment) $1$.
We remark that this final spatial scaling also plays a role in understanding the performance of
the $M/M/1$ queue as the traffic intensity $\rho$ increases toward its critical value $1$.  With appropriate scaling
of both time and space, the stochastic process $Q$ approaches reflected Brownian motion
with negative drift, while the busy period distribution has interesting behavior, in which both small values and large values play a role
\cite{AW88a,AW95,W02}.

In terms of random variables, let
\beql{b3}
Y_{\sigma} \equiv \frac{X_{\sigma/(1+\sigma)}}{1 + \sigma} =
 (1 - \rho) X_{\rho}, \quad \mbox{where} \quad \rho \equiv \frac{\sigma}{1 + \sigma} \qandq \sigma = \frac{\rho}{1 - \rho}.
\eeq
Let $b(x; \sigma)$ be the mgf of $Y_{\sigma}$.  From above,
\beql{b4}
b(x; \sigma) \equiv E[e^{x Y_{\sigma}}] = \frac{1 + 2 \sigma - x - \Psi (x)}{2 \sigma}, \quad \mbox{where} \quad  \Psi (x) \equiv \sqrt{1 - 2(1 + 2 \sigma)x + x^2}.
\eeq
Let the associated moments of
$Y_{\sigma}$ be
$m_n (\sigma) \equiv E[Y_{\sigma}^n]$, $n \ge 1$,
where $m_1 (\sigma) = 1$ for all $\sigma > 0$.  These moments, divided by $n!$, are the coefficients of the series expansion of $b(x; \sigma)$
\beql{b6}
b(x; \sigma) = \sum_{n = 0}^{\infty} \frac{m_n (\sigma) x^n}{n!}
=  \sum_{n = 0}^{\infty} b_n (\sigma) x^n, \quad \mbox{where} \quad b_n (\sigma) \equiv \frac{m_{n} (\sigma)}{n!}.
\eeq
The coefficients $b_n (\sigma)$ (and thus also
the moments $m_n (\sigma)$) are polynomials in $\sigma$ (and thus integers when $\sigma$ is an integer); the first few are:
$b_0 (\sigma) = 1$, $b_1 (\sigma) = 1$, $b_2 (\sigma) = 1 + \sigma$ and $b_3 (\sigma) = 1 + 3 \sigma + 2 \sigma^2$.

More specifically, we now relate the coefficients $b_n (\sigma)$ in \eqn{b6} to the Catalan numbers, denoted by $C_k$, (A000108)
starting with $1, 1, 2, 5, 14, 42, 132, 429, 1430$; in particular,
\beql{Catalan}
C_n \equiv \frac{1}{n+1} \binom{2n}{n} \qandq C_{n+1} = \sum_{i=0}^{n} C_i C_{n-i}, \quad n \ge 1.
\eeq
  The Catalan numbers can be
characterized in terms of their generating function
\beql{b7}
c(y) \equiv \sum_{k=0}^{\infty} C_k y^k \equiv \frac{1 - \sqrt{1 - 4y}}{2y}.
\eeq
\begin{theorem}{$($mixing moments of the busy period pdf$)$}\label{th1}  For $n \ge 1$,
\beql{b8}
b_{n+1} (\sigma) = \sum_{k = 0}^{n} \binom{n+k}{n-k} C_k \sigma^k,
\eeq
where $b_n (\sigma)$ is defined in {\em \eqn{b4}} and {\em \eqn{b6}}, and $C_k$ is the $k^{\rm th}$ Catalan number in {\em \eqn{Catalan}}.
In addition, there is a convenient recurrence relation, with
$b_{0} (\sigma) \equiv b_{1} (\sigma) \equiv 1$ and
\beql{b8a}
(n + 1) b_{n+1} (\sigma) = (2n -1)(1 + \sigma) b_n (\sigma) - (n-2) b_{n-1} (\sigma).
\eeq
\end{theorem}
We provide a proof in \S \ref{secProofs}.  For $\sigma = 1,2, 3$, the coefficient sequences are:
\begin{eqnarray}
\{b_n (1)\} & = & 1, \quad 1, \quad  2, \quad  6, \quad  22, \quad  90, \quad  394, \quad  (A155069) \nonumber \\
\{b_n (2)\} & = & 1, \quad 1, \quad  3, \quad  15, \quad  93, \quad 645, \quad 4791, \quad  (A103210) \nonumber \\
\{b_n (3)\} & = & 1, \quad 1, \quad  4, \quad  28, \quad  244, \quad  2380, \quad  24868, \quad  (A103211) \nonumber
\end{eqnarray}
Sequence $\{b_n (1)\}$ (A155069) is a relatively recent addition to the OEIS \cite{S10}, having the title ``Expansion of
$(3 - x - \sqrt{1 - 6x +x^2})/2$ in powers of $x$.''  We provide a model context.  However, sequence
$\{b_n (1)\}$ shifted one to the left, i.e., $1, 2, 6, 22, 90, \ldots$ is (A006318), corresponding to the famous
``Large Schroeder numbers.''  The ``little Schroeder numbers'' in (A001003) are obtained by dividing A006318 by $2$, i.e.,
the sequence $1, 1, 3, 11, 45, \ldots$.

The moment sequences themselves, $m_n (\sigma) \equiv n! b_n (\sigma)$ are of course
also integer sequences, but they are evidently not in the OEIS \cite{S10}.  For example, the sequences
with terms $n! b_n (1)$ and $n! b_n(1/2)$, $n \ge 0$, yielding $1,1,4,36, 528, 10800$  and $1,1,3, 18,171, 2250$, respectively, are not found.

The busy period is the first passage time of $Q$ from $1$ to $0$.  Since the first passage time pdf from each state $k$ to $0$ is the $k$-fold convolution
of the busy period pdf \cite[Theorem 3.1]{AW88a}), the moments and mixing moments of these pdf's also generate integer sequences for integer $\sigma$.

Since we find regularity by multiplying $\{b_n (\sigma)\}$ by $1/2 = \sigma/(1 + \sigma)$ when $\sigma = 1$, we are motivated to consider
the associated sequences $\{\sigma b_{n+1} (\sigma)/(1 + \sigma): n \ge 1\}$ for integers $\sigma > 1$.
\begin{corollary}\label{cor1} For $n \ge 1$,
\beq
\frac{\sigma b_{n+1} (\sigma)}{1 + \sigma} = \sum_{k = 0}^{n} (-1)^{n+k} \binom{n+k}{n-k} C_k (1 + \sigma)^k.
\eeqno
\end{corollary}

We have already observed that Corollary \ref{cor1} gives the little Schroeder numbers for $\sigma = 1$; it is also discussed in
\cite[Problem 15, p. 168]{R68}.

\subsection{The Busy-Period Stationary-Excess Distribution}\label{secExcess}

Let $B_e$ denote the busy period stationary-excess cdf associated with the busy-period cdf $B_{\rho}$, defined as in \eqn{a21}
after applying the scaling in \eqn{b3}.  Let $b_e$ be the associated stationary-excess pdf.
 We can apply Proposition \ref{propShift} to characterize the mixing moments $b_{e,n} (\sigma)$ of $b_e$.
\begin{corollary}{$($mixing moments of the busy-period stationary-excess pdf$)$}\label{cor2} The $M/M/1$ busy-period stationary excess pdf $b_e$
has mgf
\begin{eqnarray}\label{be}
b_e (x, \sigma) & \equiv &\int_{0}^{\infty} e^{x t} b_e (t) \, dt \equiv \sum_{n=0}^{\infty} b_{e,n} (\sigma)x^n = \frac{b(x; \sigma) - 1}{x} \nonumber \\
& = & \frac{2}{1 - x + \sqrt{1 - 2(1+ 2\sigma)x + x^2}} \nonumber \\
& = & \frac{c(\sigma x/(1-x)^2)}{1-x} = \frac{2}{1 - x + \Psi(x)},
\end{eqnarray}
where $c(y)$ is the generating function of the Catalan numbers in {\em \eqn{b7}} and
$\Psi$ is defined in {\em \eqn{b4}}.
For $n \ge 1$,
\beql{ex1}
b_{e,n} (\sigma) = b_{n+1} (\sigma), n\ge 1,
\eeq
 for which an expression is given in Theorem {\em \ref{th1}}.  The mean of $b_e$ is $\sigma$.
\end{corollary}

\subsection{Large Schroeder and Catalan Numbers as Moments of PDF's}

In this section we identify the pdf's whose moments are (i) the Large Schroeder numbers and (ii) the Catalan numbers.
By Corollary \ref{cor2} and our observation after Theorem \ref{th1}, the sequence $\{b_{e,n} (1)\}$ coincides with the Large Schroeder numbers (A006318).
From \cite[Theorem 5.1 and Corollary 5.2.1]{AW88a} and \cite[Theorem 4.1]{AW88b}, we can obtain the spectral representation for the pdf $b_e$.
(We need to account for the different scaling of time used there.)  That identifies the desired pdf and, by Proposition \ref{propMixing}, the mixing moments.

\begin{corollary}{$($Large Schroeder numbers as moments$)$}\label{cor3} The large Schroeder numbers arise as the moments of the mixing pdf
associated with the stationary-excess busy-period pdf $b_e$ when $\sigma = 1$.
The mixing pdf for $b_e$ as a function of $\sigma$ is
 \beq
f(y; \sigma) = \frac{\sqrt{(\tau - y)(y - \tau^{-1})}}{2 \sigma \pi y}, \quad \frac{1}{\tau} < y < \tau,
\eeqno
where $\tau \equiv 1 + 2\sigma + 2\sqrt{\sigma(1 + \sigma)}$.
\end{corollary}

We introduced the Catalan numbers before Theorem \ref{th1}.  The Catalan numbers are related
to reflected Brownian motion with
drift $-1$ and diffusion coefficient $1$, denoted by $\{R(t): t\ge0\}$.
It is the limit of the $M/M/1$ queue length process as $\rho \uparrow 1$ with appropriate scaling of
time and space \cite{AW87,W02}.  Since $E[R(t)|R(0) = 0]$ is nondecreasing,
$H_1 (t) \equiv E[R(t)|R(0) = 0]/E[R(\infty)]$, $t \ge 0$, is a cdf.
By \cite[Corollary 1.5.2]{AW87}, \cite[\S 7]{AW96} and \cite[(8.7)]{AW99a},
the pdf $h_1$ of $H_1$ has LT $\hat{h}_1 (s)$, which can be characterized as the unique fixed point of the exponential mixture operator, i.e.,
\beql{h2}
\hat{h}_1 (s) = \frac{1}{1 + s \hat{h}_1 (s)}.
\eeq
Together with the spectral representation given in \cite[Theorem 4.1]{AW88b},
that implies the following result.  Without the probability model context, the result already appears in the commentary on
(A000108) \cite{P02}.
\begin{theorem}{$($Catalan numbers as moments$)$}\label{thCat} The generating function $c(x)$ of the
Catalan numbers in {\em \eqn{Catalan}} and {\em \eqn{b7}} coincides
with the mgf of $h_1$, $\hat{h}_1 (-x)$, associated with the LT in {\em \eqn{h2}}.
As a consequence, the Catalan numbers arise as the moments of the mixing density of $h_1$,
\beql{h3}
f(y) = \frac{\sqrt{4 - y}}{2\pi \sqrt{y}}, \quad 0 < y < 4.
\eeq
\end{theorem}

If we take the Laplace transform of $f$ in \eqn{h3} and replace $s$ by $-x$, then we obtain the mgf,
which is the egf $\tilde{c} (x)$ of the Catalan numbers.
\begin{corollary}\label{corEGF}
The egf of the Catalan numbers is
\beq
\tilde{c} (x) \equiv \sum_{n=0}^{\infty} \frac{C_n x^n}{n!} = e^{2x}(I_{0} (2x) - I_{1} (2x)),
\eeqno
where $I_0$ and $I_1$ are Bessel functions.
\end{corollary}

\subsection{The $M/M/1$ Equilibrium Time to Emptiness}

From a sequence perspective, the shifting by $1$ and multiplying by $\sigma/(1+\sigma)$ following Theorem \ref{th1} and Corollary \ref{cor1}
lead us to consider the mgf
\beql{b10}
p(x; \sigma) \equiv \frac{1}{1 + \sigma} + \frac{\sigma}{1 + \sigma}\left(\frac{b(x; \sigma) - 1}{x}\right).
\eeq

The function $p(x; \sigma)$ in \eqn{b10} turns out to be the moment generating function of the equilibrium time to emptiness
in the $M/M/1$ queue, i.e., the first passage time to state $0$ by the stochastic process $Q$, assuming that $Q$ starts at time $0$ according to
its (geometric) steady-state distribution.
Hence, with probability $(1-\rho) = 1/(1 + \sigma)$ the process starts at $0$, so the first passage time is $0$.
As a consequence, there is an atom at $0$ with mass $1/(1 + \sigma)$.  Conditional on starting at a positive value, the equilibrium time to emptiness
coincides with the stationary-excess pdf $b_e$ of the busy-period pdf $b$ \cite[Theorem 3]{AW94}.  Additional characterizations appear
in \cite[Theorem 3.3]{AW88a}.  We make a connection to the exponential mixture operator in \eqn{a41}.
Paralleling \eqn{b6}, we let $p_n (\sigma)$ be the coefficient of $p(x; \sigma)$ as a gf, i.e., writing
$p (x;\sigma) \equiv \sum_{n=0}^{\infty} p_n (\sigma) x^n$.
We provide a proof in \S \ref{secProofs}.

\begin{theorem}{$($equilibrium time to emptiness$)$}\label{th2}
The mgf $p(x;\sigma)$ in {\em \eqn{b10}} can be expressed as
\beql{b15}
p(x; \sigma)  =  \frac{c((1+\sigma)x/(1+x)^2)}{1+x} = \frac{2}{1 + x + \Psi(x)},
\eeq
where $c(y)$ is the generating function of the Catalan numbers in {\em \eqn{b7}} and
$\Psi$ is defined in {\em \eqn{b4}}.  Hence,
\beql{b15a}
p_n (\sigma) = \frac{\sigma b_{n+1}(\sigma)}{1+\sigma} =  \frac{\sigma}{1 + \sigma} \sum_{k=0}^{n} \binom{n+k}{n-k} C_k \sigma^k, \quad n \ge 0,
\eeq
and the following recursion can be applied with $p_0 (\sigma) \equiv 1$ and $p_1(\sigma) \equiv \sigma$
\beql{b15b}
(n+2) p_{n+1} (\sigma) = (2n+1)(1+2\sigma) p_n (\sigma) - (n-1) p_{n-1}(\sigma), \quad n \ge 1.
\eeq
In addition, $b(x; \sigma)$ can be expressed as an exponential mixture of $p(x; \sigma)$, i.e.,
\beql{b16}
b(x; \sigma) = \frac{1}{1 - x p(x; \sigma)}.
\eeq
\end{theorem}

From \eqn{b15a} and \eqn{b15b},
 the sequence $\{p_n (\sigma): n \ge 1\}$ is an integer sequence for each positive integer $\sigma$.
For $\sigma = 1,2, 3$, the coefficient sequences are:
\begin{eqnarray}
\{p_n (1)\} & = & 1, \quad 1, \quad  3, \quad  11, \quad  45, \quad  197, \quad  (A001003) \nonumber \\
\{p_n (2)\} & = & 1, \quad 2, \quad  10, \quad  62, \quad  430, \quad 3194, \quad  (A107841) \nonumber \\
\{p_n (3)\} & = & 1, \quad 3, \quad  21, \quad  183, \quad  1785, \quad 18651, \quad   (A131763). \nonumber
\end{eqnarray}

We conclude this section by remarking that in the theory of integer sequences there also is a convenient {\em invert operator},
which can be expressed for LT's via
\begin{eqnarray}\label{invert}
\hat{g} (s) & \equiv & \sI (\hat{f}(s)) \equiv \mbox{Excess} (\mbox{exp mixture}(\hat{f} (s)) = \hat{f}_{\sE\sM, e} (s) \nonumber \\
            & = &\mbox{Excess} \left(\frac{1}{1 + s \hat{f}(s)}\right) = \frac{1}{s} \left(1 - \frac{1}{1 + s \hat{f} (s)}\right) \nonumber \\
            & = & \frac{\hat{f}(s)}{1 + s \hat{f} (s)}.
\end{eqnarray}

The following result links the three mgf's $b(x; \sigma)$, $b_e (x; \sigma)$ and $(p (x; \sigma)$; see Corollary 5.1 of \cite{AW88a}.
We also consider the excess of the excess, denoted by $b_{e.e}$ \cite{W85}.
\begin{corollary}\label{cor4}  The mgf's $b_e (x; \sigma)$, $p (x; \sigma)$ and $b(x; \sigma)$ are related by
\beq
b_e (x; \sigma) = \sI (p (x; \sigma)) = p(x; \sigma) b(x; \sigma)
\eeqno
for the invert operator in {\em \eqn{invert}} and
\beq
b_{e,e} (x; \sigma) = p_e (x; \sigma).
\eeqno
\end{corollary}

Directly from the final expression in \eqn{be}, we get a heavy-traffic limit for the scaled mgf
with our scaling in \eqn{b3}, with the limit
being the gf of the Catalan numbers.  Let $X_e (\sigma)$ be a random variable with mgf $b_e (x; \sigma)$;
i.e., $b_e (x;\sigma) \equiv E[e^{x X_e (\sigma)}]$.  Let $\Rightarrow$ denote convergence in distribution \cite{W02}.
\begin{corollary}\label{cor5}  As $\sigma \ra \infty$,
\beq
 b_e (x/\sigma; \sigma) \ra  \frac{2}{1 - \sqrt{1 - 4x}} = \frac{1 - \sqrt{1 - 4x}}{2x} \equiv c(x)
\eeqno
for $c$ in \eqn{b7}, so that
\beq
\frac{X_e (\sigma)}{\sigma} \Rightarrow X \qasq \sigma \ra \infty,
\eeqno
where $c(x) = E[e^{x X}]$.
\end{corollary}
A related heavy-traffic limit for the time-scaled busy-period stationary-excess pdf was obtained in
\cite[Theorem 1]{AW95}.  (They are consistent.)

\subsection{Continued Fractions and Hankel Transforms}

We now look at the three probability mgf's $b(x; \sigma)$, $b_e (x,; \sigma) \equiv (b(x;\sigma) - 1)/x$ and $p(x; \sigma)$ in
\eqn{b4}, \eqn{be} and \eqn{b10} from the perspective of continued fractions and Hankel transforms.
From the previous subsections, we know that the associated three random variables are closely related.
That is shown again through this new perspective.  First, the three mgf's can be represented as $S$ (Stieltjes) fractions.
(In \cite{AW99a} we found that $S$ fractions frequently arise in CF's associated with birth-and-death stochastic processes.)
The CF coefficients are given in Table \ref{tb1}.
\begin{table}[h!]
\begin{center}
\begin{tabular}{|r||llllll|}
\hline \hline
mgf             & $h_1$ & $h_2$ & $h_3$ & $h_4$ & $h_5$ & $h_6$ \\ \hline
$b(x; \sigma)$     & 1         & $\sigma$    & $1 + \sigma$ & $\sigma$     & $1+\sigma$ & $\sigma$ \\
$b_e(x; \sigma)$  & $1+ \sigma$ & $\sigma$    & $1 + \sigma$ & $\sigma$     & $1+\sigma$ & $\sigma$ \\
$p(x;\sigma)$     & $\sigma$    & $1+ \sigma$ & $\sigma$     & $1 + \sigma$ & $\sigma$   & $1+\sigma$ \\
\hline \hline
\end{tabular}
\caption{The coefficients of the continued fractions representing the three $M/M/1$ mgf's.
The pattern repeats in each row.}\label{tb1}
\end{center}
\end{table}

For the generating functions of our $M/M/1$ queueing examples, the determination of the CF
coefficients is simple because we can apply the algebraic identity
\beq
 \frac{\alpha }{1+} \frac{\gamma}{1+} \frac{\alpha}{1+} \frac{\gamma}{1+} \ldots
 = \frac{1}{2}\left(\sqrt{1 + 2(\gamma + \alpha) + (\gamma - \alpha)^2} - 1 - (\gamma - \alpha)\right)
\eeqno
for constants $\alpha$ and $\gamma$;
i.e., we can solve the equation
\beq
\zeta = \frac{\alpha}{1+} \frac{\gamma}{1+ \zeta}
\eeqno
for $\zeta$.

The Hankel transform of an integer sequence provides a useful partial characterization;
it is a many-to-one function mapping an integer sequence into another integer sequence.
(For an example of non-uniqueness, see Corollary \ref{cor6} below.)
For example, the Hankel transform of the Catalan numbers is the sequence $1,1,1, \dots$ \cite{B09, L01}.
Following \cite{CBVWJ08, DS97}, starting from
a sequence $\{\omega_n: n \ge 0 \} \equiv \omega_0, \omega_1, \omega_2, \omega_2, \dots$
with $\omega_0 \equiv 1$, let the
{\em Hankel matrix} $M^{(n)}$ be the $(n+1) \times (n+1)$ symmetric matrix with elements $M^{(n)}_{i,j} \equiv \omega_{i+j-2}$,
$0 \le i \le n$, $0 \le j \le n$.
(The first row contains the first $n+1$ elements and $M_{n+1,n+1} \equiv \omega_{2n}$.
Let $H_{2n} \equiv det(M^{(n)})$, the {\em even  Hankel determinant}.
Let the {\em Hankel transform} of the sequence $\{\omega_n: n \ge 0 \}$ above be the sequence
$\{H_{2n}: n \ge 0\}$; it starts with $H_0 = 1$.

One way to compute the Hankel transform of a sequence $\{\omega_n: n \ge 0 \}$ is to first
determine the corresponding continued fraction from the formal power series; i.e.,
starting from \eqn{power}, we determine \eqn{CF} above.  The coefficients $h_n$ appearing in \eqn{CF}
can be determined by applying the normalized Viskovatov algorithm, as given on \cite[p. 112]{CBVWJ08}.
Then we apply the iteration
\beq
H_{2n} = \left(\prod_{i = 1}^{2n} h_i \right) H_{2n - 2}, \quad n \ge 1;
\eeqno
see Theorem 1.4.10 on  of \cite[p. 23]{DS97}.  Another way to arrive at this result is via the even contraction of the CF in \eqn{CF},
as given in \cite[(12.3), p. 270]{B09}; note that $\beta_n$ in \cite{B09} is equal to $h_{2n-1} h_{2n}$, $n \ge 1$, in our notation.
Then (12.2) of \cite[(12.2)]{B09} and our iteration yield the same result.

\begin{corollary}\label{cor6}
The Hankel transforms of the sequences $\{b_n (\sigma)\}$ in Theorem {\em \ref{th1}} , $\{b_{e,n} (\sigma)\}$ in Corollary {\em \ref{cor2}}
and $\{p_n (\sigma\}$ in Theorem {\em \ref{th2}} are
\begin{eqnarray}
H_{2n} (b)   & = & \sigma^{\nu(n)}(1 + \sigma)^{\nu(n) - n}, \nonumber \\
H_{2n} (p)   & = & H_{2n} (b_e) = (\sigma + \sigma^2)^{\nu(n)}, \nonumber
\end{eqnarray}
where $\nu (n) \equiv n(n+1)/2$.
\end{corollary}

\subsection{The Stationary Waiting Time in the $M/G/1$ Queue}

The mgf $w(x; \sigma) \equiv E[e^{xW}]$ of the steady-state waiting time $W$ in the $M/G/1$ queue
can be expressed in terms of the mgf $g_e (x; \sigma)$ of the stationary-excess
of the general service-time distribution using the Pollazcek-Khintchine transform (again assuming that the mean service time is $1$) by
applying the random sum operator from \cite[Table 1]{AW96}, yielding
\beql{w1}
w(x; \sigma) = \frac{1-\rho}{1 -\rho g_e (x; \sigma)} = \frac{1}{1 + \sigma - g_e (x; \sigma)}
\eeq
\cite[\S 4.3]{R62} and \cite{AW94}.  Riordan \cite[(18a), p. 49]{R62} and Tak\'{a}cs \cite{T62} develop a nice recursion for the moments,
 which as a function of $\sigma$ becomes
\beql{w2}
E[W^{n-1} (\sigma)] = \frac{\sigma}{n} \sum_{k=2}^{n} \binom{n}{2} g_k E[W^{n-k} (\sigma)] \qandq E[W^{0} (\sigma)] \equiv 1,
\eeq
where $g_k$ is the $k^{\rm th}$ moment of the service-time cdf.  Clearly, \eqn{w2} can be the source of many integer sequences.

To illustrate,
we now consider the Catalan numbers as service-time moments, which is legitimate by Theorem \ref{thCat}.  We exploit the following result, obtained by
combining \cite[Corollaries 1.3.2 and 1.5.1]{AW87} with Theorem \ref{thCat}.
\begin{lemma} Let $h_1$ be the density of the RBM first moment function, which has the Catalan numbers as its mixing moments.  Then
\beql{h5}
\hat{h}_{1,e} {s} = \hat{h}_1 (s)^2, \quad \mbox{so that} \quad c_e(x) = c(x)^2,
\eeq
where $c(x)$ is the generating function of the Catalan numbers.
\end{lemma}

\begin{theorem}
If the service time pdf is the RBM first moment pdf $h_1$ characterized by its LT in {\em \eqn{h2}}, which has the Catalan numbers as its moments,
then the waiting time mgf in {\em \eqn{w1}} becomes
\beql{rep1}
w(x; \sigma) = \frac{1}{1 + \sigma - \sigma c(x)^2}
\eeq
and the moments of the stationary waiting time pdf are
\beql{rep2}
E[W^{n-1} (\sigma)] = \sigma \sum_{k = 2}^{n} \frac{(n-1)!}{(n-k)!} C_k E[W^{n-k} (\sigma)].
\eeq
Hence, if $\sigma$ is an integer, then $\{E[W^n (\sigma)]\}$ is an integer sequence.
\end{theorem}

For the case $\sigma = 1$, we find $\{E[W^n (1)]\} = 1, 2, 18, 252, 4776, \ldots$, which is not in the OEIS.
Riordan \cite[(19), p. 50]{R62} shows how the moments $E[W^n (\sigma)]$ can be expressed in terms of the multivariate Bell polynomials.

Let $w_n (\sigma) \equiv E[W^n (\sigma)]/n!$ be the associated mixing moments.  From \eqn{rep2}, we obtain a recursion for $w_n (\sigma)$.
\begin{corollary}\label{cor7}
A recursion for the mixing moments $w_n (\sigma)$ defined above is
\beql{rep3}
w_{n-1} (\sigma) = \sigma \sum_{k=2}^{n} C_k w_{n-k} (\sigma), \quad n \ge 1.
\eeq
\end{corollary}
For $\sigma = 1$, we get
 $\{w_n (1)\} = 1, 2, 9, 42, 199, \ldots$,
 which also is not yet in OEIS.

 \subsection{General Birth-and-Death Processes}

 The most relevant previous work connecting queueing theory to integer sequences and providing a suitable framework for generalization
 evidently is the previous work connecting birth-and-death processes to continued fractions; see \cite{AW99a, FG00} and references therein;
 continued fractions are known to be intimately connected to integer sequences.
 We have already mentioned that the stochastic process $Q \equiv \{Q(t): t\ge 0\}$, representing the number of customers in
 an $M/M/1$ queueing model at time $t$, is a birth-and-death stochastic process.  Indeed, it is a special birth-and-death process with
 constant birth rate $\lambda$ and constant death rate $\mu$.  More generally, these rates are functions of the state \cite[\S 6.3]{R07};
 $\lambda_k$ is the rate up and $\mu_k$ is the rate down when $Q(t) = k$.  The general birth-and-death process represents a more general queueing model.

 For a general birth-and-death process $X \equiv \{X(t): t \ge 0\}$, an important quantity is the probability
 \beql{return}
 P_{0,0} (t) \equiv P(X(t) = 0| X(0) = 0), \quad t \ge 0;
 \eeq
 see \cite[\S 11]{AW88a}, \cite[\S 3.1]{FG00} and \cite{K72}.
From \cite[(4.1), p. 771]{FG00}, the LT of $P_{0,0} (t)$, $\hat{P}_{0,0} (s)$, has a
representation as an $S$ fraction
\begin{eqnarray}
\hat{P}_{0,0} (s) & =  &  \frac{1}{s+} \frac{\lambda_0}{1+} \frac{\mu_1} {s+} \frac{\lambda_1}{1+} \frac{\mu_2}{s+} \ldots \nonumber \\
                  & =  &  \frac{z}{1+} \frac{\lambda_0 z}{1+} \frac{\mu_1 z} {1+} \frac{\lambda_1 z}{1+} \frac{\mu_2 z}{1+} \ldots
\end{eqnarray}
where $z \equiv 1/s$.
(The last line follows from \cite[(1.5)]{AW99a} using the sequence $\{c_n: n \ge 0\}$ with $c_{2n} \equiv 1$ and $c_{2n+1} \equiv z$, $n \ge 1$.)
We can find the corresponding power series
\beq
\hat{P}_{0,0} (s) = z(1 - p_1 z + p_2 z^2 - p_3 z^3 + \ldots );
\eeqno
see \cite[(3.9), (3.10), p. 397]{AW99a}.  The sequence $\{1, p_1, p_2, p_3, \ldots\}$ may be an integer sequence.

For the special case of the $M/M/1$ queue, with our scaling we have $\lambda_i = \sigma$ and $\mu_{i+1} = 1 + \sigma$ for all $i$, $i \ge 0$,
so that
\beq
\hat{P}_{0,0} (s) = z \left(\frac{1}{1+} \frac{\sigma z}{1+} \frac{(1 + \sigma) z} {1+} \frac{\sigma z}{1+} \frac{(1 + \sigma) z}{1+} \ldots \right).
\eeqno
Then from Table \ref{tb1} we have the relation
\beql{return2}
\hat{P}_{0,0} (s) = z \hat{p} (z),
\eeq
where $\hat{p} (s)$ is the Laplace transform of the equilibrium time to emptiness.

\section{Proofs}\label{secProofs}

In this concluding section we provide proofs for Theorems \ref{th1} and \ref{th2}.

\subsection{Proof of Theorem \ref{th1}.}

We give two proofs.  The first is direct; the second applies \cite{B09a}.

\paragraph{Direct proof of \eqn{b8}.}

From \cite[pp. 232--233]{T67}, after a change of scale and notation, we have
\beq
b_{n+1} = \left(\frac{1 + \sigma}{n+1}\right) \sum_{k=1}^{n} \binom{n+k}{k}\binom{n-1}{k-1} \sigma^{k-1}.
\eeqno
(This equation is also given by \cite[(21), p. 151]{R68} after we make the identification
that $H_n (z) = p_n (\sigma) = \sigma b_{n+1} (\sigma)/(1 + \sigma)$.)  On the right side,
move $(1+ \sigma)$ inside the sum and identify the coefficients of $\sigma^k$, obtaining
\beq
\left(\frac{1}{n+1}\right) \left(\binom{n+k+1}{k+1}\binom{n-1}{k} + \binom{n+k}{k}\binom{n-1}{k-1}\right)
= \binom{n+k}{n-k}\binom{2k}{k}\left(\frac{1}{k+1}\right).~~~\bsq
\eeqno

\paragraph{Direct proof of \eqn{b8a}.}

We follow the technique of \cite[p. 107]{R62} to establish a three-term recursion.  By the successive differentiation
with respect to $x$ in \eqn{b4}, we find that
\beq
\Psi(x)^2 b''(x; \sigma) = (1+\sigma)(1 + 2 \sigma - x)b' (x; \sigma) + (1+\sigma)^2 b(x; \sigma).
\eeqno
In this equation, make the following substitutions:
\begin{eqnarray}
b(x; \sigma) = \sum_{n=0}^{\infty} b_n (\sigma) x^n, \nonumber \\
b'(x; \sigma) = \sum_{n=0}^{\infty} (n+1)b_{n+1} (\sigma) x^n, \nonumber \\
b''(x; \sigma) = \sum_{n=0}^{\infty} (n+1)(n+2) b_{n+2} (\sigma) x^n. \nonumber
\end{eqnarray}
Then collect and equate the coefficients of $x^{n-1}$ and the result follows.~~~\bsq

\paragraph{Application of \cite{B09a}.}

The idea in this second proof is to directly make connection to the Catalan numbers via their generating function $c$ in \eqn{b7}.
Starting from \eqn{b4}, we apply \eqn{be}, \eqn{b10} and the proof of \eqn{b15} in Theorem \ref{th2} below to conclude that
\beql{k1}
b_e (x; \sigma) = \frac{1}{1-x}c\left(\frac{\sigma x}{(1 - x)^2}\right) \qandq p (x; \sigma) = \frac{1}{1+x}c\left(\frac{(1 +\sigma) x}{(1 + x)^2}\right).
\eeq
We then apply \cite[Proposition 15, p.12]{B09a} (with change of variables $y \equiv \sigma x$ and $y \equiv (1+\sigma) x$, respectively) to immediately deduce the conclusions
\beql{k2}
b_{e, n} (\sigma) = \sum_{k=0}^{n} \binom{n+k}{n-k} C_k \sigma^k \qandq p_{n} (\sigma) = \sum_{k=0}^{n} \binom{n+k}{n-k} (-1)^{n+k} C_k (1 +\sigma)^k
\eeq
given in Corollary \ref{cor2} and Theorem \ref{th2}.  We then can apply \eqn{ex1} and \eqn{be} to get the conclusion in Theorem \ref{th1} for $b$ from $b_e$.
By this line of reasoning, we establish Theorems \ref{th1} and \ref{th2} and Corollaries \ref{cor1} and \ref{cor2} all at once,
and have a new perspective on their relationship.

\subsection{Proof of Theorem \ref{th2}.}

\paragraph{Proof of \eqn{b15}.}

From the representation of $p(x; \sigma)$ in \eqn{b10}, we have
\begin{eqnarray}
p(x; \sigma) & = &\frac{1}{1+\sigma} + \frac{\sigma}{1+\sigma}\left(\frac{2}{1 - x + \Psi(x)}\right) \nonumber \\
& = & \frac{1}{1+\sigma} + \frac{\sigma}{1+\sigma}\left(\frac{1 - x - \Psi(x)}{2 \sigma x}\right) \nonumber \\
& = & \frac{1 + x - \Psi (x)}{2 (1 + \sigma) x} = \frac{2}{1 + x + \Psi(x)}. \nonumber
\end{eqnarray}
On the other hand, by \eqn{b7},
\begin{eqnarray}
\frac{c((1+\sigma)x/(1 +x)^2)}{1+x} & = &\frac{2}{1 + x + (1 + x)\sqrt{1 - 4(1+\sigma) x/(1+x)^2}}  \nonumber \\
& = & \frac{2}{1 + x + \Psi(x)}.~~~\bsq \nonumber
\end{eqnarray}

\paragraph{Proof of \eqn{b15a}.}

Combine \eqn{b8} and \eqn{b10}.~~~\bsq

\paragraph{Proof of \eqn{b15b}.}

Combine \eqn{b8a} and \eqn{b10}.~~~\bsq

\paragraph{Proof of \eqn{b16}.}

From \eqn{b4},
\begin{eqnarray}
b(x; \sigma) & = & \frac{1 + 2 \sigma -x - \Psi(x)}{2 \sigma} = \frac{2(1 + \sigma)}{1 + 2\sigma - x + \Psi(x)} \nonumber \\
&=& \frac{1}{1 - x\left([1+x - \Psi(x)]/[2(1 + \sigma)x]\right)} = \frac{1}{1 - x p(x; \sigma)}.~~~\bsq \nonumber
\end{eqnarray}

\section{Acknowledgment}

The second author was supported by NSF Grant CMMI 0948190.

\pagebreak

\bigskip
\hrule
\bigskip

\noindent 2000 {\it Mathematics Subject Classification}:
Primary 11B25, 60K25.

\noindent \emph{Keywords: }  integer sequence, moment sequence,
probability, probability density function, probability moment
sequence, the moment problem, queue, queueing, Laplace
transforms, moment generating functions, generating functions,
continued fractions, Hankel transform, renewal theory, birth-and-death process,
stationary excess distribution, stationary lifetime
distribution, exponential mixture, single-server queue, $M/M/1$
queue, $M/G/1$ queue, busy period, equilibrium time to
emptiness, waiting time.

\bigskip
\hrule
\bigskip

\noindent (Concerned with sequences
\seqnum{A000108}, \seqnum{A001003}, \seqnum{A006318}, \seqnum{A093526}, \seqnum{A103210}, \seqnum{A103211},
\seqnum{A107841}, \seqnum{A122525},  \seqnum{A131763}, \seqnum{A137216}, \seqnum{A155069},
.)

\bigskip
\hrule
\bigskip

\vspace*{+.1in} \noindent Received February 19, 2010; Revision
May 3, 2010.


\begin{thebibliography}{5}

\bibitem{ACLW95}
J. Abate, G. L. Choudhury, D. M. Lucantoni and W. Whitt, Asymptotic analysis of tail probabilities based on the computation of moments,
{\em Ann. Appl. Prob.} {\bf 5} (1995), 983--1007.



\bibitem{AW87}
J. Abate and W. Whitt, Transient behavior of regulated Brownian motion,
{\em Adv. Appl. Prob.} {\bf 19} (1987), 560--598.


\bibitem{AW88a}
J. Abate and W. Whitt, Transient behavior of the $M/M/1$ queue via Laplace transforms,
{\em Adv. Appl. Prob.} {\bf 20} (1988), 145--178.

\bibitem{AW88b}
J. Abate and W. Whitt, Simple spectral representations for the $M/M/1$ queue,
{\em Queueing Systems} {\bf 3} (1988), 321--346.

\bibitem{AW94}
J. Abate and W. Whitt, Transient behavior of the $M/G/1$ workload process,
{\em Operations Research} {\bf 42} (1994), 750--764.

\bibitem{AW95}
J. Abate and W. Whitt, Limits and approximations for the busy-period distribution in single-server queues,
{\em Prob. Engr. Inform. Sci.} {\bf 9} (1995), 581--602.


\bibitem{AW96}
J. Abate and W. Whitt, An operational calculus for probability distributions via Laplace transforms,
{\em Adv. Appl. Prob.} {\bf 28} (1996), 75--113.


\bibitem{AW99a}
J. Abate and W. Whitt, Computing Laplace transforms for numerical inversion via continued fractions,
{\em INFORMS Journal on Computing} {\bf 11} (1999), 394--405.

\bibitem{AW99b}
J. Abate and W. Whitt, Explicit $M/G/1$ waiting-time distributions for a class of long-tail service-time distributions..
{\em Operations Research Letters} {\bf 25} (1999), 25--31.

\bibitem{A65}
N. I. Akhiezer, {\em The Classical Moment Problem and Some Related Questions in Analysis}, Oliver and Boyd, London, 1965.

\bibitem{A04}
M. H. Albert, R. E. L. Aldred, M. D. Atkinson, H. P. van Ditmarsch, C. C. Handley
and D. A. Holton, Restricted permutations and queue junping, {\em Discrete Math.} {\bf 287} (2004), 129--133.

\bibitem{AB94}
M. D. Atkinson, R. Beals, Priority queues and permutations,  {\em SIAM J. Comput.} {\bf 23} (1994), 1225--1230.

\bibitem{B09}
P. Barry, {\em A Study of Integer Sequences, Riordan Arrays, Pascal-like Arrays and Hankel Transforms},
PhD dissertation Department of Mathematics, University College Cork, Ireland, 2009.

\bibitem{B09a} P. Barry,  Continued fractions and
    transformations of integer sequences,
    {\em J. Integer Sequences} {\bf 12} (2009), 09.7.6.


\bibitem{CBVWJ08}
A. Cuyt, V. Brevik Petersen, B. Verdonk, H. Waadeland and W. B. Jones,
{\em Handbook of Continued fractions for Special Functions}, Springer, 2008.

\bibitem{CL96}
G. L. Choudhury and D. M. Lucantoni, Numerical computation of the moments of a probability distribution
from its transform, {\em Oper. Res.} {\bf 44}  (1996), 368--381.

\bibitem{DS97}
H. Dette and W. J. Studden, {\em The Theory of Canonical Moments with Applications to Statistics, Probability and Analysis},
Wiley, New York, 1997.

\bibitem{FG00}
P. Flajolet and F. Guillemin, The formal theory of birth-and-death processes, lattice path combinatorics and continued fractions,
{\em Adv. Appl. Prob.} {\bf 32} (2000), 750--778.



\bibitem{K72}
J. F. C. Kingman, {\em Regenerative Phonomena},
Wiley, New York, 1972.

\bibitem{L01}
J. W. Layman, The Hankel transform and some of its properties,
{\em J. Integer Sequences} {\bf 4} (2001), 01.1.5.


\bibitem{P02}
K. A. Penson, An integral representation of $C_n$,
commentary on (A000108), OEIS, 2002.

\bibitem{R58}
J. Riordan, {\em An Introduction to Combinatorial Analysis},
 Wiley, New York, 1958.

 \bibitem{R62}
J. Riordan, {\em Stochastic Service Systems},
 Wiley, New York, 1962.

 \bibitem{R68}
J. Riordan, {\em Combinatorial Identities},
 Wiley, New York, 1968.

  \bibitem{R07}
S. M. Ross, {\em Introduction to Probability Models}, ninth edition,
Academic Press, New York, 2007.

 \bibitem{S10}
N. J. A. Sloane,
\href{http://www.research.att.com/~njas/sequences/}{The On-Line
Encylopedia of Integer Sequences} (OEIS), 2010.



 \bibitem{S05}
R. P. Stanley,
Queue problems revisited, Suomen Tehtavaniekat (Proceedings of the Finnish Chess Problem Society) {\bf 59} (2005), 193--203.


 \bibitem{T62}
L. Tak\'{a}cs, {\em Introduction to the Theory of Queues}, Oxford University Press, New York, 1962.

 \bibitem{T67}
L. Tak\'{a}cs,
{\em Combinatorial Methods in the Theory of Stochastic Processes},  Wiley, New York, 1967.

 \bibitem{T90}
 L. Tak\'{a}cs, On the number of distinct forests,  {\em SIAM J. Discrete Math.} {\bf 3} (1990), 574--581.

\bibitem{W85}
W. Whitt,
 The renewal-process stationary-excess operator. {\em J. Appl. Prob.} {\bf 22} (1985), 156--167.

\bibitem{W02}
W. Whitt, {\em Stochastic-Process Limits}, Springer, New York, 2002.


\end{thebibliography}
\end{document}